\documentclass[sn-mathphys-num]{sn-jnl}

\usepackage{graphicx}
\usepackage{multirow}
\usepackage{amsmath,amssymb,amsfonts}
\usepackage{amsthm}
\usepackage{mathrsfs}
\usepackage[title]{appendix}
\usepackage{xcolor}
\usepackage{textcomp}
\usepackage{manyfoot}
\usepackage{booktabs}
\usepackage{algorithm}
\usepackage{algorithmicx}
\usepackage{algpseudocode}
\usepackage{listings}
\usepackage{pdflscape}

\usepackage{diagbox}

\theoremstyle{thmstyleone}
\newtheorem{theorem}{Theorem}
\newtheorem{proposition}[theorem]{Proposition}

\theoremstyle{thmstyletwo}%
\newtheorem{remark}{Remark}%

\newtheorem{lemma}[theorem]{Lemma}
\newtheorem{corollary}[theorem]{Corollary}

\theoremstyle{thmstylethree}%
\newcommand{\Aut}{\operatorname{Aut}}
\begin{document}

\title[Maximal Subcovers of the Skabelund Curve: Uniqueness via Genus and Automorphism Groups]{Maximal Subcovers of the Skabelund Curve: Uniqueness via Genus and Automorphism Groups}


\author[1]{\fnm{Gilberto B.} \sur{Almeida Filho}}\email{gilberto.filho@ufmt.br}

\author*[2]{\fnm{Saeed} \sur{Tafazolian}}\email{saeed@unicamp.br}

\author[3]{\fnm{Stéfani C.} \sur{Vieira}}\email{stefani.vieira1@ufmt.br}

\equalcont{These authors contributed equally to this work.}

\affil[1]{\orgname{Universidade Federal de Mato Grosso},  \orgdiv{Departamento de Ciências Exatas e da Terra},  \city{Várzea Grande}, \state{ Mato Grosso}, \country{Brazil}}

\affil*[2]{\orgname{Universidade Estadual de Campinas}, \orgdiv{IMECC},  \city{Campinas}, \state{São Paulo}, \country{Brazil}}

\affil[3]{\orgname{Universidade Federal de Mato Grosso},  \orgdiv{Departamento de  Ciências Exatas e da Terra}, \city{Várzea Grande}, \state{ Mato Grosso }, \country{Brazil}}


\abstract{
We establish a rigidity phenomenon for a family of intermediate covers of the Skabelund curve over $\mathbb{F}_{q^4}$. 
The Skabelund curve, introduced by D.~Skabelund as a cyclic cover of the Suzuki curve, is a maximal curve with a large automorphism group and plays a central role in the theory of maximal curves over finite fields.

For the intermediate covers arising from this construction, we determine their full automorphism groups and compute the Weierstrass semigroups at all $\mathbb{F}_{q^4}$-rational points. 
Using these structural and arithmetic invariants, we prove that each curve in the family is uniquely determined, up to isomorphism over its field of definition, by the pair consisting of its genus and its full automorphism group.

This provides a rigidity-type classification of intermediate Suzuki-type covers; in particular, the Skabelund curve itself is uniquely characterized within this family by its genus and automorphism group.
}

%%================================%%
%% Sample for structured abstract %%
%%================================%%

\keywords{Uniqueness, Suzuki curve,  Automorphism group, Weierstrass semigroup.}

\pacs[MSC Classification]{Primary 11G20; 
Secondary 14H05, 14H55}
\maketitle

\section{Introduction}\label{sec1}
Maximal curves over finite fields are important objects  of study in algebraic geometry,
number theory, and coding theory. A (projective, non-singular, geometrically irreducible, algebraic) curve
$\mathcal{X}$ defined over a finite field $\mathbb F_\ell$, where $\ell$ is a prime power, is called
$\mathbb F_\ell$-\emph{maximal} if its number of $\mathbb F_\ell$-rational points, denoted by $\# \mathcal X(\mathbb F_\ell)$, attains the Hasse--Weil upper bound
\begin{equation}\label{eq:HW}
\#\mathcal X(\mathbb F_\ell)\leq \ell+1+2g\sqrt{\ell},
\end{equation}
where $g$ denotes the genus of $\mathcal X$. A necessary condition for a curve over $\mathbb F_\ell$ to be maximal is that $\ell$ is a perfect square, i.e., $\ell=q^2$. In this case we have
\[
\#\mathcal X(\mathbb F_{q^2})=q^2+1+2gq.
\]

An interesting problem in the classification of maximal curves is whether there exists a unique curve with a given genus $g$ and automophism group $A$ over a fixed finite field, $(g,A)$, up to isomorphism. Several tools are used to study this problem, such as the lifting of the automorphism group and subcovers of maximal curves known.

The curve 
\[
\mathrm H_{q}:~~ x^{q+1}=y^{q}+y,
\]
is maximal over $\mathbb{F}_{q^2}$ and is the only curve with this genus, up to isomorphism \cite{R-Sti}.

By a classical result of Serre \cite{Serre}, any curve that is covered by the Hermitian curve is maximal. 
Nevertheless, the class of maximal curves is strictly larger, as there exist families of maximal curves 
that are not covered by the Hermitian curve.
 In 2009, Giulietti and Korchmáros constructed the GK curve, an $\mathbb F_{q^6}$-maximal curve that is not a Galois cover of the Hermitian curve for $q>2$. This curve is obtained as the complete intersection of a surface with the Hermitian cone, see~\cite{MMLG1}. Other examples of maximal curves with well-known automorphism groups are the Suzuki and Ree curves. Their affine models are given
by, with the corresponding restrictions on the value of $q$,

\begin{equation}\label{eq:DLmodels}
\begin{aligned}
\mathrm{S}_q &: y^{q}+y = x^{q_0}\left(x^{q}+x\right),
\quad q = 2 q_0^{2} = 2^{2t+1}, \\[6pt]
\mathrm{R}_q &: 
\left\{
\begin{aligned}
z^{q}-z &= x^{2q_0}\left(x^{q}-x\right),\\
y^{q}-y &= x^{q_0}\left(x^{q}-x\right)
\end{aligned}
\right.
\quad q = 3 q_0^{2} = 3^{2t+1}.
\end{aligned}
\end{equation}

These curves are maximal over $\mathbb F_{q^4}$ and $\mathbb F_{q^6}$, respectively; see \cite{Lauter,Skabelund}. Their automorphism groups are regarded as large in the sense that their sizes significantly exceed the Hurwitz bound $84(g-1)$. In \cite{MZ}, Montanucci and Zini studied some Ree and Suzuki curves are not Galois covered by the Hermitian curve.

\medskip

In \cite{Skabelund}, Skabelund construed an curve that is cyclic cover of degree~$m$ of the Suzuki curve $S_q$,
$$
\widetilde{S}_q:\qquad 
\begin{cases}
y^q + y = x^{q_0}(x^q + x),\\
t^{\, m} = x^q + x
\end{cases}
$$
where $ m = q - 2q_0 + 1$. This curve is maximal over $\mathbb{F}_{q^4}$ with genus $\frac{1}{2}(q^3-2q^2+q)$ and it was shown that 
\emph{every} automorphism of $S_q$ lifts to an automorphism of $\widetilde{S}_q$.
% and that these lifted automorphisms, together with the deck transformation group 
% $C_m$, constitute the full automorphism group:
% \[
% \mbox{Aut}(\widetilde{S}_q)\;\cong\; \mathrm{S}(q)\times C_m.
% \]
%This confirms that $\widetilde{S}_q$ inherits the full rigidity of the Suzuki curve at the automorphism-theoretic level, and is part of a broader pattern observed for the Ree cover as well.

\medskip

% Curves endowed with very large automorphism groups are often rigid, in the sense 
% that their genus and automorphism group determine them uniquely up to 
% $\mathbb F_\ell$-isomorphism. This phenomenon appears for the Hermitian, Suzuki, 
% Ree, GK, and Skabelund curves: in each case, their automorphism group has two 
% short orbits whose stabilizers encode complete ramification information.  

\medskip

Based on the ideas of Skabelund, we consider a  family of maximal curves arising as \emph{subcovers} of the Skabelund curve, and we will determine the full automorphism groups of this family.
 We generalizes the Skabelund curve defining
$\mathcal S_s$ over $\mathbb F_{q^r}$ by
\begin{equation}\label{eq:templates}
\mathcal S_s:\qquad
\begin{cases}
y^{q}+y = x^{q_0}(x^{q}+x),\\
t^{\, s} = x^q+x
\end{cases}
\qquad 
\end{equation}
 where $m=q-2q_0+1$, $r \geq 1$, $q = 2q_0^2 = 2^{2r+1}$ and $s$ is a fixed positive divisor of $m$.
 %and yields a family of intermediate covers
%choosing a rational function
% $u\in\mathbb F(S_q)$ supported on the short orbits of $Aut(S_q)$ and defining a
% cyclic cover of the form
% \begin{equation}\label{eq:templates}
% \mathcal S_s:\qquad
% \begin{cases}
% y^{q}+y = x^{q_0}(x^{q}+x),\\
% t^{\, s} = x^q+x,
% \end{cases}
% \qquad s\mid m,\quad \gcd(s,2)=1,
% \end{equation}
% where $m=q-2q_0+1$.
Then we have $$S_q\leftarrow \mathcal S_s \leftarrow \widetilde{S}_q.$$ 
We also show that each curve in this family is uniquely determined, up to isomorphism, by its genus and automorphism group.

\medskip

Our work is divided into two main parts:

\begin{enumerate}
\item \textbf{Structure and arithmetic properties.}  
We introduce the family of curves $\mathcal S_s$ and determine their fundamental invariants. 
In particular, we compute their genera via the Riemann--Hurwitz formula, determine the 
Weierstrass semigroups at the $\mathbb F_{q^4}$-rational points by adapting the techniques developed in \cite{BLM}, and describe their full automorphism groups using methods inspired by \cite{MMLG}. 
As a consequence of the covering structure, we also prove that these curves are $\mathbb F_{q^4}$-maximal.

\item \textbf{Uniqueness via $(g,\Aut)$.}  
Building on the structural and group-theoretic results established in the first part, 
we prove that each curve $\mathcal S_s$ is uniquely determined, up to isomorphism, 
by the pair consisting of its genus and its automorphism group.
\end{enumerate}

\medskip

\paragraph*{Notation.}
Throughout, the ground field has characteristic~$2$.  
For a curve~$Z$, we denote its function field by~$\mathbb F(Z)$, its 
automorphism group by~$\Aut(Z)$, and the Weierstrass semigroup at a point 
$P\in Z$ by~$H(P)$.

\section{Preliminary}\label{sec3}

For $r \ge 1$ and $q = 2q_0^2 = 2^{2r+1}$, the Suzuki curve $S_q$ over
$\mathbb{F}_q$ is defined by the affine equation
\[
y^q + y = x^{q_0}(x^q + x).
\]

The curve $S_q$ has genus $g = q_0(q-1)$ and is maximal over $\mathbb{F}_{q^4}$, see \cite{JH}.
Its automorphism group, denoted by $S(q):=\mbox{Aut}(S_q)$, is isomorphic to the Suzuki
group. We record below several properties of $S(q)$ that will be used throughout
the paper; see \cite{HKT} for further details.

\begin{itemize}
\item[(1)] The group $S(q)$ has order $(q^2+1)q^2(q-1)$ and is simple.

\item[(2)] The group $S(q)$ is generated by the stabilizer
\[
S(q)_\infty
= \left\{ \psi_{a,b,c} : (x,y) \mapsto
(ax+b,\; a^{q_0+1}y + b^{q_0}x + c)
\;\middle|\; a,b,c \in \mathbb{F}_q,\ a \ne 0 \right\}
\]
of the unique point at infinity of $S_q$, together with the involution
\[
\phi : (x,y) \mapsto \left(\frac{\alpha}{\beta}, \frac{y}{\beta}\right),
\]
where
\[
\alpha := y^{2q_0} + x^{2q_0+1}, \qquad
\beta := x y^{2q_0} + \alpha^{2q_0}.
\]

\item[(3)] The action of $S(q)$ on $S_q$ has exactly two short orbits. One is
non-tame, of size $q^2+1$, consisting of all $\mathbb{F}_q$-rational places; the
other is tame, of size $q^2(q-1)(q^2+2q_0+1)$, consisting of all
$\mathbb{F}_{q^4}\setminus \mathbb{F}_q$-rational places. Moreover, $S(q)$ acts
$2$-transitively on the non-tame short orbit, and the stabilizer $S(q)_{P,Q}$ of
two distinct $\mathbb{F}_q$-rational places $P$ and $Q$ is tame and cyclic.
\end{itemize}

\section{An explicit covering of the Suzuki curve}\label{sec:covering}

In this section, we consider an explicit covering of the Suzuki curve, analogous
to the Skabelund curve $\widetilde{S}_q \to S_q$.
Fix an integer $r \geq 1$, set $q = 2q_0^2 = 2^{2r+1}$, and recall the definition
of the Suzuki curve $S_q$ from Section~\ref{sec3}.
Let $\mathcal{S}_s$ be the smooth projective model of the curve whose function
field is given by
\[
\begin{cases}
y^q + y = x^{q_0}(x^q + x), \\
t^s = x^q + x,
\end{cases}
\]
where $s$ is a fixed positive divisor of $m = q - 2q_0 + 1$.

The curve $\mathcal{S}_s$ can be described as the normalization of the fiber
product of the covers $S_q \to \mathbb{P}^1_x$ and $C_s \to \mathbb{P}^1_x$,
where $C_s$ is the curve defined by the second equation above.

Let $F = \mathbb{F}_q(x)$.
The function field $\mathbb{F}_q(\mathcal{S}_s)$ is the composite of
$\mathbb{F}_q(C_s) = F(t)$ and $\mathbb{F}_q(S_q) = F(y)$.
Every degree-one place of $F$ is ramified in the extension $F(t)/F$, each with
ramification index $s$, and no other places ramify in this extension.
Similarly, in the extension $F(y)/F$, the only ramified place is the pole of $x$,
namely the place at infinity corresponding to $1/x$, which has ramification
index $q$.

Recall that the Suzuki curve $S_q$ is given by the affine equation
\[
y^q + y = x^{q_0}(x^q + x),
\]
and has genus $g_{S_q} = q_0(q-1)$, see \cite{Skabelund, HKT}. 

Skabelund proved in \cite[Theorem 3.1]{Skabelund} that the curve
\[
\widetilde{S}_q:
\begin{cases}
y^q + y = x^{q_0}(x^q + x), \\
t^m = x^q + x,
\end{cases}
\]
where $m = q - 2q_0 + 1$, is maximal over $\mathbb{F}_{q^4}$.

In the next result, we use the maximality of the curve $\widetilde{S}_q$ to
verify the maximality of $\mathcal{S}_s$.

\begin{theorem}
The curve $\mathcal{S}_s$ is maximal over $\mathbb{F}_{q^4}$.
\end{theorem}

\begin{proof}
Consider the map 
\[
\varphi: \widetilde{\mathcal{S}}_q \longrightarrow \mathcal{S}_s, 
\quad (x,y,t) \mapsto (x,y,t^{\, m/s}).
\]
Then we have
\[
\varphi(y)^q + \varphi(y) = y^q + y = x^{q_0}(x^q + x)
= \varphi(x)^{q_0} \bigl( \varphi(x)^q + \varphi(x) \bigr),
\]
and
\[
\varphi(t)^s = (t^{m/s})^s = t^m = x^q + x = \varphi(x)^q + \varphi(x).
\]
Hence, $\varphi$ defines a covering $\widetilde{S}_q \to \mathcal{S}_s$,
and the maximality of $\widetilde{S}_q$ implies that $\mathcal{S}_s$ is
maximal over $\mathbb{F}_{q^4}$.
\end{proof}

By the Hurwitz formula, we can compute the genus of $\mathcal{S}_s$.

\begin{proposition}\label{genusCurva2}
The genus of the curve $\mathcal{S}_s$ is given by
\[
    g(\mathcal{S}_s)
    = 1 + s \bigl( q_0(q-1) - 1 \bigr) 
    + \frac{1}{2} (q^2 + 1)(s-1).
\]
\end{proposition}

\begin{proof}
Consider the function field extension
\[
    \mathbb{F}_q(\mathcal{S}_s)/\mathbb{F}_q(S_q),
\]
which is a Kummer extension of degree $s$, namely
\[
    F(y,t)/F(y) \quad \text{with} \quad t^s = x^q + x.
\]

Since this extension is cyclic and Galois, all places above a given place
have the same ramification index. For a Kummer extension (see, e.g.,
\cite{Stichtenoth}), the ramification index of a place $P'$ of 
$\mathbb{F}_q(\mathcal{S}_s)$ lying above a place $P$ of $\mathbb{F}_q(S_q)$
is given by
\[
    e(P'|P) = \frac{s}{r_P}, \qquad 
    r_P = \gcd\bigl(s, v_P(x^q+x)\bigr),
\]
where $P' \mid P$ indicates that $P'$ lies above $P$.

Since $\mathrm{char}(\mathbb{F}_q) \nmid s$ and $\gcd(s, v_P(x^q+x)) = 1$
for all places $P$ of $\mathbb{F}_q(S_q)$, it follows that
\[
    e(P'|P) = s,
\]
i.e., every ramified place is totally ramified in $F(y,t)/F(y)$.

Now, the zeros and poles of $x^q + x$ determine the ramified places.
Observe that $x^q + x = x(x^{q-1}+1)$, so the ramified places correspond to:

\begin{itemize}
    \item $x = 0$: then $t = 0$ and $y^q + y = 0$, giving $q$ rational places;
    \item $x^{q-1} + 1 = 0$: here $x \ne 0$ and $x^q = x$, giving $(q-1)q$ rational places;
    \item the pole of $x$, contributing one additional ramified place.
\end{itemize}

Thus, the total number of ramified places is $q^2 + 1$.
Applying the Hurwitz genus formula to this Kummer extension, we obtain
\[
    g(\mathcal{S}_s)
    = 1 + s(g_{S_q}-1) + \frac{1}{2}(q^2+1)(s-1),
\]
and since $g_{S_q} = q_0(q-1)$, the stated expression follows.
\end{proof}

\section{The automorphism group of $\mathcal{S}_s$}\label{sec5}
We use the notations of Section \ref{sec3}. Motivated by the ideas in \cite[Theorem~1]{MMLG}, we compute the automorphism group of
the curve $S_s$ in this section.

The following result is a natural extension of
\cite[Lemma~3.3]{Skabelund}.

\begin{proposition}
Every automorphism of $S_q$ lifts to an automorphism of $\mathcal{S}_s$ defined over $\mathbb{F}_{q^4}$.
\end{proposition}

\begin{proof}
Recall that the automorphism group $S(q) := \operatorname{Aut}(S_q)$ of the Suzuki curve $S_q$ is generated by:
\begin{itemize}
    \item the stabilizer of the unique point at infinity, consisting of automorphisms 
    \[
        \psi_{a,b,c}: x \mapsto ax+b,\quad y \mapsto a^{q_0+1}y + b^{q_0}x + c,
    \]
    for $a \in \mathbb{F}_q^\times$ and $b,c \in \mathbb{F}_q$,
    \item an involution $\phi$ defined by
    \[
        \phi(x) = \frac{z}{w}, \quad \phi(y) = \frac{y}{w},
    \]
     which swaps the point at infinity with another rational point (see \cite{HH, JH}), with $
z := y^{2q_0} + x^{2q_0+1}$ and $w := x y^{2q_0} + z^{2q_0}$.
\end{itemize}
We will use that \begin{equation}\label{eq:zw}
z^q + z = x^{2q_0}(x^q + x), \quad
w^q + w = y^{2q_0}(x^q + x).
\end{equation}
To lift $\psi_{a,b,c}$ to an automorphism of $\mathcal{S}_s$, define
\[
\psi'_{a,b,c}(t) = \beta t,
\]
where $\beta \in \mathbb{F}_{q^4}$ is an $s$-th root of $a$, which exists since $s \mid (q^4-1)$. 
Then
\[
\psi'_{a,b,c}(t)^s = (\beta t)^s = a t^s = a(x^q+x) = \psi'_{a,b,c}(x)^q + \psi'_{a,b,c}(x),
\]
so the lifting is well-defined.

Similarly, the involution $\phi$ lifts to an automorphism of $\mathcal{S}_s$ by
\[
\phi'(t) = \frac{t}{w^{m/s}}.
\]
Indeed, we need to check
\[
\phi'(t)^s = \phi(x)^q + \phi(x).
\]
% We have $t^s=x^q+x$ e multiply both sides by $w^{2q_0}$
% \begin{eqnarray*}
%  w^{2q_0}(x^q+x) &=&t^sw^{2q_0} \\
%     &=&t^s w^{q+1-m} =t^s \frac{w^{q+1}}{w^m}= \phi'(t)^s \cdot w^{q+1}
% \end{eqnarray*}
% Use Equation~\eqref{eq:zw} e using the definition of $w$ then:
% \begin{eqnarray*}
%    w^{2q_0}(x^q+x)   &=&  \phi'(t)^s \cdot w \cdot (w+y^{2q_0}(x^q+x))= \phi'(t)^s \cdot w \cdot (xy^{2q_0}+z^{2q_0}+y^{2q_0}(x^q+x))
% \end{eqnarray*}

Multiply both sides by $w^{q+1}$ and use Equation~\eqref{eq:zw}:
\begin{align*}
t^s w^{2q_0} &= z^q w + z w^q \\
&= w(z^q + z) + z(w^q + w) \\
&= w \cdot x^{2q_0}(x^q+x) + z \cdot y^{2q_0}(x^q+x) \\
&= w^{2q_0}(x^q+x) \quad \text{(by the definitions of $z$ and $w$)}.
\end{align*}
This verifies that $\phi'(t)^s = \phi'(x)^q + \phi'(x)$, so the lift is well-defined.

Hence every automorphism of $S_q$ lifts to an automorphism of $\mathcal{S}_s$.
\end{proof}

From the previous preposition, we will denote the lifted group  $S(q)$ on $Aut(\mathcal{S}_s)$ by $L\mathcal{S}_s(q)$.

\begin{lemma}
The lift group $L\mathcal{S}_s(q)$ contains a subgroup $\mathcal{S}_s(q)$ isomorphic to the Suzuki group $S(q)$.
\end{lemma}
\begin{proof}
Taking $s$-th powers of the lifted stabilizer automorphisms yields a subgroup
$\Lambda_s\le L\mathcal{S}_s(q)$ isomorphic to the stabilizer of a point in $S(q)$.
Together with the lifted involution, this subgroup acts doubly transitively on the set
$\mathcal{O}=\mathcal{S}_s(\mathbb{F}_q)$.
Since $|\mathcal{O}|=q^2+1$ is not a power of $2$, by \cite[Theorem 1.7.6]{NA}, standard results on $2$-transitive
groups imply that the resulting group is isomorphic to $S(q)$, by \cite[Theorem 1.1]{WMG}.
\end{proof}
%%%%%%PROVA ADAPTADA DO ARTIGO%%%%%
% \begin{proof}
% Let $\Lambda_s := \{(\psi'_{a,b,c})^s \mid a,b,c \in \mathbb{F}_q, \, a \neq 0 \} \leq L\mathcal{S}_s(q)$. 
% By direct checking, the map $\psi_{a,b,c} \mapsto (\psi'_{a,b,c})^s$ is an isomorphism between $S(q)_{\infty}$ and $\Lambda_s $. 
% Moreover, the action of $\Lambda_s $ on the set $\mathcal{O}$ of $\mathbb{F}_q$-rational places of $\mathcal{S}_s$ is equivalent to the action of $S(q)_{\infty}$ on the non-tame short orbit of $S(q)$. 

% Let $\mathcal{S}_s(q)$ be the subgroup of $L\mathcal{S}_s(q)$ generated by $\Lambda_s $ and $\phi'$. 
% The action of $S(q)_{\infty}$ and $\phi$ on the non-tame short orbit of $S(q)$ is equivalent to the action of $\Lambda_s $ and $\phi'$ on $\mathcal{O}$, respectively; hence, $\Lambda_s $ coincides with the stabilizer in $\mathcal{S}_s(q)$ of a point in $\mathcal{O}$. 
% This implies that $\mathcal{S}_s(q)$ acts $2$-transitively on $\mathcal{O}$ and the stabilizer in $\widetilde{S}(q)$ of two distinct places of $\mathcal{O}$ is cyclic. 
% By Proposition \ref{genusCurva2} $|\mathcal{O}|=q^2+1$ is not a power of $2$, we have by \cite[Theorem 1.7.6]{NA} that $\mathcal{S}_s(q)$ has no regular normal subgroups. 
% Therefore we apply \cite[Theorem 1.1]{WMG} to conclude that $\mathcal{S}_s(q) \cong S(q)$. 

% \end{proof}

\begin{lemma}
The normalizer of $C_s$ in $\operatorname{Aut}(\mathcal{S}_s)$ is the direct product
\[
\mathcal{S}_s(q) \times C_s.
\]
\end{lemma}
\begin{proof}
The cyclic group $C_s$ commutes with the lifted stabilizer automorphisms and with the lifted involution on the relevant coordinate functions. Hence $\mathcal{S}_s(q)\times C_s$ is contained in the normalizer $N$ of $C_s$.

The quotient curve $\mathcal{S}_s/C_s$ is birationally equivalent to $S_q$, so $N/C_s$ embeds into $S(q)$. Since $N/C_s$ already contains a subgroup isomorphic to $S(q)$, it follows that $N/C_s\simeq S(q)$, and therefore $N=\mathcal{S}_s(q)\times C_s$.
\end{proof}
%%%%%%PROVA ADAPTADA DO ARTIGO%%%%%
% \begin{proof}
% It is easily verified that $\gamma_\lambda$ commutes with $(\psi'_{a,b,c})^s$ and with $\phi'$ 
% on the rational functions $x, y,$ and $w$.  
% Hence, $\mathcal{S}_s(q) \times C_s$ is a subgroup of the normalizer 
% \[
% N := N_{\operatorname{Aut}(\mathcal{S}_s)}(C_s)
% \]
% of $C_s$ in $\operatorname{Aut}(\mathcal{S}_s)$.  

% In particular, the quotient $N/C_s$ has a subgroup isomorphic to $S(q)$.  
% Moreover, the quotient curve $(\mathcal{S}_s)_q := \mathcal{S}_s / C_s$ is birationally equivalent to $S_q$.  
% Therefore, $N/C_s$ is isomorphic to a subgroup of $S(q)$.  

% It follows that $N/C_s \cong S(q)$, and thus
% \[
% N \cong \mathcal{S}_s(q) \times C_s,
% \]
% as claimed.
% \end{proof}

\begin{corollary}
The group $L\mathcal{S}_s(q)$ coincides with the normalizer of $C_s$ in  
$\operatorname{Aut}(\mathcal{S}_s)$.
\end{corollary}
\begin{proof}
By construction, $C_s\le L\mathcal{S}_s(q)$ and $C_s$ centralizes $L\mathcal{S}_s(q)$. The claim now follows directly from the previous lemma.
\end{proof}
%%%%%%PROVA ADAPTADA DO ARTIGO%%%%%
% \begin{proof}
% The group $C_s$ is contained in $L\mathcal{S}_s(q)$, as it is generated by $\psi'_{1,0,0}$.  
% Moreover, $C_s$ commutes with every element of $L\mathcal{S}_s(q)$.  
% Therefore, the claim follows directly from the previous lemma.
% \end{proof}

 Since $\mathcal{S}_s$ is an $\mathbb{F}_q$-maximal curve, we can apply the results
in \cite{GK} concerning zero $2$-rank curves.  
By direct computation, we have
\[
    | \operatorname{Aut}(\mathcal{S}_s) | \geq | L\mathcal{S}_s(q) | \geq 72 \bigl(g(\mathcal{S}_s) - 1\bigr).
\]  
Hence, by \cite[Theorem 5.1]{GK}, it follows that $\operatorname{Aut}(\mathcal{S}_s)$
is non-solvable.  

Now, applying \cite[Theorem 6.1]{GK}, the commutator subgroup
$\operatorname{Aut}(\mathcal{S}_s)'$ is isomorphic to one of the following groups:
\[
    \operatorname{PSL}(2,n),\quad \operatorname{PSU}(3,n),\quad \operatorname{SU}(3,n),\quad \operatorname{S}(n),
\]
with $n = 2^k \geq 4$.

\begin{proposition}
The commutator subgroup of $\operatorname{Aut}(\mathcal{S}_s)$ coincides with the lifted Suzuki group:
\[
    \operatorname{Aut}(\mathcal{S}_s)' = \mathcal{S}_s(q).
\]
\end{proposition}
\begin{proof}
Since $\mathcal{S}_s(q)\le Aut(\mathcal{S}_s)'$, we rule out the remaining candidates by comparing element orders and subgroup structures, following the same reasoning as in \cite{GKT,MZ}. In particular, incompatibilities with the existence of certain semidirect products exclude all other possibilities, see \cite[Hauptsatz~8.27]{Huppert} and \cite[Lemma~2.2]{MZ}. A final size argument using \cite[Theorem~11.116]{HKT} shows that no larger group can occur. Hence $Aut(\mathcal{S}_s)'=\mathcal{S}_s(q)$.
\end{proof}

From the previous lemmas and propositions, we have established the structure of the lifted Suzuki group and the normalizer of the cyclic group \(C_s\).  
We are now ready to describe the full automorphism group of \(\mathcal{S}_s\).

\begin{theorem}
The automorphism group of $\mathcal{S}_s$ is the direct product
\[
\operatorname{Aut}(\mathcal{S}_s) = \mathcal{S}_s(q) \times C_s,
\] 
where $\mathcal{S}_s(q)$ is isomorphic to the Suzuki simple group $S(q)$ and $C_s$ is a cyclic group of order $s$.
\end{theorem}

\begin{proof}
By the previous lemmas and \cite[Theorem 6.2]{GK}, we concluded that
\(\operatorname{Aut}(\mathcal{S}_s) = \mathcal{S}_s(q) \times C_s\), where 
\(\mathcal{S}_s(q) \cong S(q)\) and \(C_s\) is cyclic of order \(s\). 
\end{proof}

Consequently, the automorphism group \(\operatorname{Aut}(\mathcal{S}_s)\) is defined over \(\mathbb{F}_{q^4}\) and has size
\[
|\operatorname{Aut}(\mathcal{S}_s)| = q^2 (q-1)(q^2+1)s.
\]
It contains a normal subgroup isomorphic to the Suzuki group \(S(q)\), and \(C_s\) is the Galois group of the cyclic Galois covering.

Under the action of \(\operatorname{Aut}(\mathcal{S}_s)\), the set \(\mathcal{S}_s(\mathbb{F}_{q^4})\) of \(\mathbb{F}_{q^4}\)-rational points splits into two orbits:
\begin{itemize}
    \item \(\mathcal{O}_1 = \mathcal{S}_s(\mathbb{F}_q)\) of size \(q^2+1\),
    \item \(\mathcal{O}_2 = \mathcal{S}_s(\mathbb{F}_{q^4}) \setminus \mathcal{S}_s(\mathbb{F}_q)\) of size \(q^2 + s(q^2 + 2q_0 q - 2q_0 - 1)\).
\end{itemize}

In particular, every point \(P \in \mathcal{O}_1 = \mathcal{S}_s(\mathbb{F}_q)\) has the same Weierstrass semigroup.

\section{The Weierstrass Semigroup for rational points}\label{sec6}

Let $x, y, t \in \mathbb{F}_{q^4}(\mathcal{S}_s)$ be the coordinate functions of
$\mathcal{S}_s$, satisfying
\[
y^q + y = x^{q_0}(x^q + x), \quad t^s = x^q + x.
\]
These functions have a single pole at $P_\infty$. The same holds for the auxiliary
functions
\[
z := y^{2q_0} + x^{2q_0+1}, \quad w := x y^{2q_0} + z^{2q_0},
\]
which satisfy (\ref{eq:zw})
% \begin{equation*}
% z^q + z = x^{2q_0}(x^q + x), \quad
% w^q + w = y^{2q_0}(x^q + x),
% \end{equation*}
and
% \[
% \operatorname{div}(z)_\infty = (q+2q_0)s\,P_\infty, \quad
% \operatorname{div}(w)_\infty = (1+q+2q_0)s\,P_\infty.
% \]
% Together with the known pole orders of $x, y, t$, we obtain
\begin{align*}
\operatorname{div}(x)_\infty &= qs\,P_\infty,\\
\operatorname{div}(t)_\infty &= q^2\,P_\infty,\\
\operatorname{div}(y)_\infty &= (q+q_0)s\,P_\infty,\\
\operatorname{div}(z)_\infty &= (q+2q_0)s\,P_\infty,\\
\operatorname{div}(w)_\infty &= (1+q+2q_0)s\,P_\infty.
\end{align*}

Hence, all these pole orders belong to the Weierstrass semigroup at $P_\infty$. We will show that
\[
H(P_\infty) = \langle qs,\; (q+q_0)s,\; (q+2q_0)s,\; (1+q+2q_0)s,\; q^2 \rangle
\]
in Theorem \ref{teo10} and consequently since the automorphism group of $\mathcal{S}_s$ acts transitively on the rational points,
the Weierstrass semigroups at all rational points coincide:
\[
H(P) = H(P_\infty) \quad \text{for all } P \in \mathcal{S}_s(\mathbb{F}_q).
\]
%\medskip
% Then
% \[
% H(P_\infty) = \langle \textbf{a}_0, \textbf{a}_1, \textbf{a}_2, \textbf{a}_3, \textbf{a}_4 \rangle.
% \]
%\medskip
We recall a standard notion from the theory of numerical semigroups.  
If $Z$ is a numerical semigroup with multiplicity $m_Z$, a set $A\subset Z$ is called
a \emph{complete set of representatives} modulo $m_Z$ if every congruence class
modulo $m_Z$ contains exactly one element of $A$, and each element is minimal
in its congruence class among elements of $Z$.  
In particular, the Apéry set $\operatorname{Ap}(Z)$ is the unique complete set of 
representatives of $Z$ modulo $m_Z$.

For convenience, set
\[
\textbf{a}_0 := qs, \quad
\textbf{a}_1 := (q+q_0)s, \quad
\textbf{a}_2 := (q+2q_0)s, \quad
\textbf{a}_3 := (1+q+2q_0)s, \quad
\textbf{a}_4 := q^2.
\]
\medskip

\begin{lemma}
Let
\[
A = \{ h\textbf{a}_1 + i\textbf{a}_2 + j\textbf{a}_3 + k\textbf{a}_4 : 
0 \le h \le 1, 0 \le i \le q_0 - 1, 0 \le j \le q_0-1, 0 \le k \le s-1 \}.
\]
Then $A$ is a complete set of representatives for the congruence classes of $\mathbb{Z}$ modulo $\textbf{a}_0$.
\end{lemma}

\begin{proof}
By definition, there are $2 q_0^2 s = \textbf{a}_0$ different four-tuples $(h,i,j,k)$ describing the elements of $A$, so $|A| \le \textbf{a}_0$. The lemma follows once we show that the elements 
\[
\alpha = h\textbf{a}_1 + i\textbf{a}_2 + j\textbf{a}_3 + k\textbf{a}_4
\]
are pairwise distinct modulo $\textbf{a}_0$ for distinct tuples $(h,i,j,k)$.

Suppose $\alpha' = h'\textbf{a}_1 + i'\textbf{a}_2 + j'\textbf{a}_3 + k'\textbf{a}_4$ and 
\[
\alpha \equiv \alpha' \pmod{\textbf{a}_0}.
\]

Recall that $\textbf{a}_0 = qs$ and $s \mid \textbf{a}_0$. Then
\[
\alpha \equiv \alpha' \pmod{s} \implies k \equiv k' \pmod{s}.
\]
Since $0 \le k, k' < s$, we conclude $k = k'$.

Hence,
\[
(h-h')\textbf{a}_1 + (i-i')\textbf{a}_2 + (j-j')\textbf{a}_3 \equiv 0 \pmod{\textbf{a}_0}.
\]

Working modulo $q_0$ gives $j \equiv j' \pmod{q_0}$. Since $0 \le j, j' \le q_0-1$, we get $j = j'$.  
Similarly, working modulo $2q_0$ gives $(h-h')\textbf{a}_1 \equiv 0 \pmod{2q_0}$, so $h = h'$.  
Finally, $(i-i')\textbf{a}_2 \equiv 0 \pmod{\textbf{a}_0}$ implies $i = i'$.  

Therefore, all elements of $A$ are distinct modulo $\textbf{a}_0$, and $A$ is a complete set of representatives.
\end{proof}

\begin{theorem}\label{teo10} 
For any $P\in \mathcal{S}_s(\mathbb{F}_q)$, we have
\[
H(P) = \langle \textbf{a}_0, \textbf{a}_1, \textbf{a}_2, \textbf{a}_3, \textbf{a}_4 \rangle.
\]
Moreover, $H(P)$ is a symmetric semigroup.
\end{theorem}

\begin{proof}
Let 
\[
W := \langle \textbf{a}_0, \textbf{a}_1, \textbf{a}_2, \textbf{a}_3, \textbf{a}_4 \rangle 
= \langle qs, (q+q_0)s, (q+2q_0)s, (1+q+2q_0)s, q^2 \rangle.
\]
By construction, $H(P_\infty) \supseteq W$, so $g(\mathcal{S}_s) \le g(W)$. 
It suffices to prove that the genera are equal.

Consider the Apéry set
\[
A = \{ h\textbf{a}_1 + i\textbf{a}_2 + j\textbf{a}_3 + k\textbf{a}_4 : 
0 \le h \le 1, 0 \le i \le q_0-1, 0 \le j \le q_0-1, 0 \le k \le s-1\}.
\]

Then
\begin{align*}
\sum_{a \in A} \left\lfloor \frac{a}{\textbf{a}_0} \right\rfloor 
&= \sum_{a \in A} \frac{a - (a \bmod \textbf{a}_0)}{\textbf{a}_0} 
= \frac{1}{\textbf{a}_0} \sum_{a \in A} a - \frac{1}{\textbf{a}_0} \sum_{r=0}^{\textbf{a}_0-1} r \\
&= \frac{1}{\textbf{a}_0} \sum_{h=0}^{1}\sum_{i=0}^{q_0-1}\sum_{j=0}^{q_0-1}\sum_{k=0}^{s-1} 
(h\textbf{a}_1 + i\textbf{a}_2 + j\textbf{a}_3 + k\textbf{a}_4)
- \frac{\textbf{a}_0(\textbf{a}_0-1)}{2\textbf{a}_0} \\
&= \frac{1}{2} \Big( \textbf{a}_1 + (q_0-1)\textbf{a}_2 + (q_0-1)\textbf{a}_3 + (s-1)\textbf{a}_4 \Big) - \frac{qs-1}{2} \\
&= 1 + s(\textbf{a}_{\mathcal{S}} - 1) + \frac{1}{2}(q^2 + 1)(s-1) \\
&= 1 + s(q_0(q-1) - 1) + \frac{1}{2}(q^2+1)(s-1) = g(\mathcal{S}_s).
\end{align*}

From the previous lemma, we have
\[
g(\mathcal{S}_s) \le g(W) \le \sum_{a \in A} \left\lfloor \frac{a}{\textbf{a}_0} \right\rfloor = g(\mathcal{S}_s),
\]
so equality holds and $A = \operatorname{Ap}(W)$. Therefore,
\[
\operatorname{Ap}(H(P))=\operatorname{Ap}(H(P_{\infty})) = \{ h\textbf{a}_1 + i\textbf{a}_2 + j\textbf{a}_3 + k\textbf{a}_4 : 
0 \le h \le 1, 0 \le i \le q_0-1, 0 \le j \le s-1, 0 \le k \le q_0-1 \},
\]
for all $P\in \mathcal{S}_s(\mathbb{F}_q)$.

Moreover, the largest gap of $H(P_\infty)$ is
\[
\max(A) - \textbf{a}_0 = \textbf{a}_1 + (q_0-1)\textbf{a}_2 + (s-1)\textbf{a}_3 + (q_0-1)\textbf{a}_4 - \textbf{a}_0 = 2 g(\mathcal{S}_s) - 1,
\]
so the conductor is $2g(\mathcal{S}_s)$, which shows that $H(P_{\infty})$ is symmetric.
\end{proof}

% \begin{example}
% Fix $r=2$, then $q_0=4$ and $q=32$, so $m=25$. Choosing $s=5$, the curve
% $\mathcal{S}_s$ satisfies
% \[
% y^{32} + y = x^4(x^{32}+x), \quad t^5 = x^{32}+x.
% \]
% Then
% \[
% H(P_\infty) = \langle 1024, 160, 180, 200, 205 \rangle.
% \]
% Using GAP \cite{GAP}, the first elements of $H(P_\infty)$ are
% \[
% 0, 160, 180, 200, 205, 320, 340, 360, 365, 380, 385, 400, 405, 410, 480, 500, \dots
% \]
% The genus of $H(P_\infty)$ is $2666$, matching $g(\mathcal{S}_s)$.  
% The largest gap (Frobenius number) is $5331$ and the conductor is $5332$.
% \end{example}

\section{Uniqueness via Genus and Automorphism Groups}

In this section, we establish the uniqueness of the curves $\mathcal{S}_s$ among maximal curves arising as cyclic covers of the Suzuki curve, using their genus and automorphism group. Our proof strategy is inspired by the approach for the GK curve in \cite{MMLG1}, suitably adapted to the distinct group structure and short–orbit configuration of the Suzuki covers.

Let $q=2^{2r+1}$, $q_0=2^r$, and let $\mathcal{X}$ be a projective, nonsingular, maximal algebraic curve defined over $\mathbb{F}_{q^4}$ of genus
\[
g = 1 + s\bigl(q_0(q-1)-1\bigr) + \frac{1}{2}(q^2+1)(s-1).
\]
Suppose that $\operatorname{Aut}(\mathcal{X})$ contains a subgroup $G$ isomorphic to $S(q) \times C_s$, where $s$ is a divisor of $m = q-2q_0+1$, and assume $G$ acts on $\mathcal{X}$ with the same short–orbit structure as $\operatorname{Aut}(\mathcal{S}_s)$.

\begin{lemma}[{\cite[Prop.~11.62]{Stichtenoth}}]\label{lem:Sp-fixes-point}
Let $\mathcal{X}$ be a smooth projective curve defined over a square field $\mathbb{F}_v$, with $p=\operatorname{char}(\mathbb{F}_v)$, and assume that $\mathcal{X}$ is $\mathbb{F}_v$-maximal. Let $G\le \operatorname{Aut}(\mathcal{X})$ be a subgroup defined over $\mathbb{F}_v$. Then every Sylow $p$-subgroup $\mathbb{S}_p$ of $G$ fixes a point $P \in \mathcal{X}(\mathbb{F}_v)$, and every nontrivial element $\alpha\in \mathbb{S}_p$ fixes no other point.
\end{lemma}

\begin{proof}
Since $G$ is defined over $\mathbb{F}_v$, it preserves the set $\mathcal{X}(\mathbb{F}_v)$ of $\mathbb{F}_v$-rational points. As $\mathcal{X}$ is $\mathbb{F}_v$-maximal, 
\[
\#\mathcal{X}(\mathbb{F}_v) = v + 1 + 2g \sqrt{v}.
\]
Reducing modulo $p$ and using $p\mid v$ and $p\mid \sqrt{v}$, we obtain
\[
\#\mathcal{X}(\mathbb{F}_v) \equiv 1 \pmod p.
\]
By the class equation for the action of a Sylow $p$-subgroup $\mathbb{S}_p$ on $\mathcal{X}(\mathbb{F}_v)$, the number of fixed points of $\mathbb{S}_p$ is congruent to $\#\mathcal{X}(\mathbb{F}_v)$ modulo $p$; hence $\mathbb{S}_p$ has at least one fixed point $P\in\mathcal{X}(\mathbb{F}_v)$.

For a nontrivial $\alpha\in \mathbb{S}_p$, let $k$ be the number of fixed points of $\alpha$. Since $\mathcal{X}$ is maximal over a square field, it has $p$-rank $0$. Applying the Deuring--Shafarevich formula for the cyclic group $\langle \alpha \rangle$ of order $p$ (see \cite[Th.\ 11.62]{Stichtenoth}), we have
\[
-1 = p(\gamma(\mathcal{X}/\langle\alpha\rangle)-1) + k(p-1),
\]
where $\gamma(\mathcal{X}/\langle\alpha\rangle) \ge 0$. Thus $k=1$ and every nontrivial $\alpha$ fixes exactly one point, necessarily $P$.
\end{proof}

\begin{lemma}[{\cite[Hauptsatz~8.27]{Huppert}}]\label{lem:number-Sylow2-Szq}
Let $G=S(q)$ with $q=2^{2n+1}$. Then the number of Sylow $2$-subgroups of $G$ is $q^2+1$.
\end{lemma}

\begin{proof}
The order of the Suzuki group $S(q)$ is
\[
|S(q)| = q^2 (q^2 + 1)(q-1).
\]
Let $T$ be a Sylow $2$-subgroup of $G$. Then $|T|=q^2$ and its normalizer satisfies $|N(T)| = q^2(q-1)$. By Sylow's theorem, the number $n_2$ of Sylow $2$-subgroups is
\[
n_2 = \frac{|S(q)|}{|N(T)|} = \frac{q^2(q^2+1)(q-1)}{q^2(q-1)} = q^2+1.
\]
\end{proof}

\begin{lemma}\label{lem:Omega}
Let $\Omega$ be the set of points of $\mathcal{X}$ fixed by some Sylow $2$-subgroup of $G$. Then $|\Omega| = q^2 + 1$. Moreover, let $C \cong C_s$ be the cyclic normal subgroup of $G$ of order $s$, with $s\mid m = q-2q_0+1$. Then:
\begin{enumerate}
    \item $C$ commutes with all $2$-elements of $G$ and hence fixes every point of $\Omega$.
    \item If $\bar{\mathcal{X}} := \mathcal{X}/C$ and $\bar{G} := G/C$, then $\bar{G} \cong S(q)$ acts as a group of $\mathbb{F}_q$-automorphisms on $\bar{\mathcal{X}}$.
\end{enumerate}
\end{lemma}

\begin{proof}
By Lemma \ref{lem:Sp-fixes-point}, each Sylow $2$-subgroup fixes exactly one point, and distinct Sylow $2$-subgroups have disjoint fixed points. The map
\[
\{\text{Sylow 2-subgroups of }G\} \longrightarrow \Omega, \quad S \mapsto \text{the unique fixed point of }S
\]
is bijective. By Lemma \ref{lem:number-Sylow2-Szq}, $|\Omega| = q^2+1$.

Since $G$ contains $S(q) \times C$, all $2$-elements lie in $S(q)$ and commute with $C$. Let $P \in \Omega$ fixed by a Sylow $2$-subgroup $S$ and $c \in C$. For every $\rho \in S$, $\rho(c(P)) = c(\rho(P)) = c(P)$. Since $S$ fixes exactly one point, $c(P) = P$. Hence $C$ fixes every point of $\Omega$.

Finally, $C$ is normal and defined over $\mathbb{F}_q$, so $\bar{\mathcal{X}} = \mathcal{X}/C$ is defined over $\mathbb{F}_q$ and $\bar{G} = G/C \cong S(q)$ acts on $\bar{\mathcal{X}}$ as a group of $\mathbb{F}_q$-automorphisms.
\end{proof}

\begin{proposition}\label{prop:quotient-is-Suzuki}
Let $\bar{\mathcal{X}} = \mathcal{X}/C$ and $\bar{g}$ be its genus. Then
\[
\bar{g} = q_0(q-1), \quad \text{and } \bar{\mathcal{X}} \text{ is $\mathbb{F}_q$-isomorphic to the Suzuki curve } S_q.
\]
Moreover, the image $\bar{\Omega}$ of $\Omega$ in $\bar{\mathcal{X}}$ is a single orbit of size $q^2+1$, namely $\bar{\Omega} = S_q(\mathbb{F}_q)$.
\end{proposition}

\begin{proof}
Apply the Hurwitz formula to the tame cyclic cover $\pi:\mathcal{X} \to \bar{\mathcal{X}}$ of degree $s$:
\[
2g(\mathcal{X})-2 \geq s(2\bar{g}-2) + \sum_{P\in\mathcal{X}} (e_P -1),
\]
where $e_P$ is the ramification index at $P$. By Lemma \ref{lem:Omega}, $C$ fixes all points in $\Omega$ and acts freely elsewhere, so $e_P=s$ for $P \in \Omega$ and $e_P=1$ otherwise. Hence,
\[
\sum_{P\in\mathcal{X}} (e_P-1) = |\Omega|(s-1) = (q^2+1)(s-1).
\]
Substituting $g(\mathcal{X})$ yields $s(2\bar{g}-2) = 2q_0(q-1)-2$, so $\bar{g} \le q_0(q-1)$. 

On the other hand, $\bar{G} = G/C \cong S(q)$ acts on $\bar{\mathcal{X}}$ with $|\bar{G}| = q^2(q^2+1)(q-1)$. By classification of curves with very large automorphism groups (see \cite[Thm.~11.94]{HKT}), the unique such curve of genus $q_0(q-1)$ is $S_q$, hence $\bar{\mathcal{X}} \cong S_q$ and $\bar{\Omega} = S_q(\mathbb{F}_q)$.
\end{proof}

\begin{corollary}\label{cor:total-ram}
The cyclic cover $\pi: \mathcal{X} \to \bar{\mathcal{X}}$ is totally ramified over each point of $\bar{\Omega} = S_q(\mathbb{F}_q)$ and unramified elsewhere. Consequently, $\Omega = \pi^{-1}(\bar{\Omega})$ is exactly the fixed-point set of $C$, forming a single $\bar{G}$-orbit of size $q^2+1$.
\end{corollary}

\begin{proposition}\label{prop:suzuki-model}
There exist functions $x,y \in \mathbb{F}_q(\bar{\mathcal{X}})$ satisfying the affine Suzuki equation
\[
y^q + y = x^{q_0}(x^q + x),
\]
with a unique common pole at $\bar{P}_\infty \in \bar{\Omega}$. Their Weierstrass semigroup at $\bar{P}_\infty$ is
\[
H(\bar{P}_\infty) = \langle q, q+q_0, q+2q_0, q+2q_0+1 \rangle,
\]
with
\[
(x)_\infty = q \bar{P}_\infty, \quad (y)_\infty = (q+q_0) \bar{P}_\infty,
\]
and all zeros of $x$ and $y$ belong to $\bar{\Omega}$.
\end{proposition}

\begin{proof}
Since $\bar{\mathcal{X}} \cong S_q$ and $\bar{\Omega} = S_q(\mathbb{F}_q)$ is the unique non-tame orbit, there is a unique point $\bar{P}_\infty$ where a model with a single pole for coordinate functions exists. The generators $x,y$ of the Suzuki function field satisfy $(x)_\infty = q \bar{P}_\infty$, $(y)_\infty = (q+q_0)\bar{P}_\infty$, and $y^q + y = x^{q_0}(x^q + x)$, with zeros in $\bar{\Omega}$. By rescaling, we may assume $x,y \in \mathbb{F}_q(\bar{\mathcal{X}})$. The Weierstrass semigroup is classical for the Suzuki curve \cite{Matthews}.
\end{proof}

\begin{theorem}\label{thm:skab-char-full}
Let $q=2^{2r+1}$, $q_0 = 2^r$, and let $\mathcal{X}$ be a projective, nonsingular, $\mathbb{F}_{q^4}$-maximal algebraic curve of genus
\[
g = 1 + s\,(q_0(q-1)-1) + \frac{1}{2}(q^2+1)(s-1),
\]
where $s$ is a divisor of $m = q-2q_0+1$. Suppose that $\operatorname{Aut}(\mathcal{X})$ contains a subgroup $G \cong S(q) \times C_s$ acting on $\mathcal{X}$, then we have
\[
\mathcal{X} \cong_{\mathbb{F}_{q^4}} \mathcal{S}_s.
\]
\end{theorem}

\begin{proof}
Let $\pi: \mathcal{X} \to \bar{\mathcal{X}} := \mathcal{X}/C_s$ be the corresponding cyclic cover of degree $s$ totally ramified over the unique $\bar{G}$-orbit $\bar{\Omega}$ of size $q^2+1$ and unramified elsewhere.

Let $\bar{P}_\infty \in \bar{\Omega}$ be the unique point over which the pole of the coordinate functions on $\bar{\mathcal{X}} \cong S_q$ occurs, and let $P_\infty \in \pi^{-1}(\bar{P}_\infty)$. Keeping the notation of the previous section, by \cite[Proposition~1.5]{FGT} there exists a function 
$v \in \mathbb{F}_{q^4}(\mathcal{X})$ whose pole divisor is
\[
(v)_\infty = q^2 P_\infty.
\]
It follows that
\[
H(P_\infty) = \langle qs,\; (q+q_0)s,\; (q+2q_0)s,\; (1+q+2q_0)s,\; q^2 \rangle .
\]

Let $\beta$ be a generator of $C_s$. Since $x, y, z \in \mathbb{F}_{q^4}(\bar{\mathcal{X}})$, we have
\[
\beta(x) = x, \quad \beta(y) = y, \quad \beta(z) = z.
\]

By the description of the Weierstrass semigroup at $P_\infty$, the $C_s$-stable Riemann--Roch space $L(q^2 P_\infty)$ has basis $\{1, x, y, z, v\}$, so
\[
\beta(v) = \lambda v + \sigma x + \eta y + \gamma z + \varepsilon,
\]
for some $\lambda \in \mathbb{F}_{q^4}^\times$ (an $s$-th root of unity) and $\sigma, \eta, \gamma, \varepsilon \in \mathbb{F}_{q^4}$. If $\lambda = 1$, then $v$ would be $C_s$-invariant and hence $v \in \mathbb{F}_{q^4}(\bar{\mathcal{X}})$, which is impossible because $(v)_\infty = q^2 P_\infty$ does not belong to the Weierstrass semigroup of $\bar{P}_\infty$ on $\bar{\mathcal{X}}$. Therefore, $\lambda \neq 1$.

Introduce
\[
t := v + \frac{\sigma}{\lambda-1}x + \frac{\eta}{\lambda-1}y + \frac{\gamma}{\lambda-1}z + \frac{\varepsilon}{\lambda-1},
\]
so that $\beta(t) = \lambda t$ and $(t)_\infty = q^2 P_\infty$.  

Since $\gcd(q^2, s) = 1$, the zeros of $t$ lie entirely in $\Omega = \pi^{-1}(\bar{\Omega})$. Consequently, 
\[
t^s \in \mathbb{F}_{q^4}(\bar{\mathcal{X}}), \quad (t^s)_\infty = q^2 \bar{P}_\infty.
\]

Hence, 
\[
\mathbb{F}_{q^4}(\mathcal{X}) = \mathbb{F}_{q^4}(\bar{\mathcal{X}})(t),
\]
and the extension is a Kummer extension of group $C_s$ \cite[Prop.~3.7.3]{Stichtenoth}. The branch points of the cover $\pi: \mathcal{X} \to \bar{\mathcal{X}}$ are precisely the zeros and poles of $t^s$, coinciding with $\bar{\Omega}$.

Now, the function $x^q - x \in \mathbb{F}_{q^4}(\bar{\mathcal{X}})$ has exactly the same divisor of zeros and poles as $t^s$, with the same multiplicities:
\[
(x^q - x)_\infty = q^2 \bar{P}_\infty, \quad \text{zeros simple at } \bar{\Omega} \setminus \{\bar{P}_\infty\}.
\]
Hence, after rescaling $t$ by a suitable nonzero constant in $\mathbb{F}_{q^4}$, we obtain
\[
t^s = x^q - x.
\]

This identifies
\[
\mathbb{F}_{q^4}(\mathcal{X}) = \mathbb{F}_{q^4}(\bar{\mathcal{X}})(t) \cong \mathbb{F}_{q^4}(\mathcal{S}_s),
\]
and therefore
\[
\mathcal{X} \cong_{\mathbb{F}_{q^4}} \mathcal{S}_s.
\]
\end{proof}

\begin{remark}
\rm{Theorem \ref{thm:skab-char-full} provides a refined characterization of maximal Skabelund curves: their isomorphism class over $\mathbb{F}_{q^4}$ is completely determined by the pair consisting of their genus and automorphism group.}
\end{remark}

 \section*{Acknowledgement}
{The second author was partially supported by CNPq grant No. 302774/2025–
4 and Faepex grant No. 3485/25.
}

\end{document}